\documentclass{article}
\usepackage{amsmath}
\usepackage{amsfonts}
\usepackage{graphicx}

\newcommand{\RR}{\mathbb{R}}

\renewcommand{\a}{\alpha}

\newtheorem{theorem}{Theorem}

\newtheorem{definition}[theorem]{Definition}

\newtheorem{lemma}[theorem]{Lemma}

\newtheorem{proposition}[theorem]{Proposition}
\newtheorem{remark}[theorem]{Remark}

\begin{document}

\title{ON THE LIENARD'S TYPE EQUATION: AN ICON OF THE NONLINEAR ANALYSIS}
\author{
\and \textbf{Juan E. N\'apoles Vald\'es}\\UNNE, FaCENA\\ Ave. Libertad 5450, Corrientes 3400, Argentina\\\texttt{jnapoles@exa.unne.edu.ar}\\
and\\
UTN-FRRE, French 414, Resistencia, Chaco 3500, Argentina \\\texttt{jnapoles@frre.utn.edu.ar}\\
}

\maketitle

\begin{abstract}
In this note, we review the latest qualitative results, referring to the Li\'enard Equation, in the framework of non-conformable, generalized and fractional differential operators.
\end{abstract}



\section{Preliminars}

The Li\'enard equation is named after the French physicist Alfred-Marie Li\'enard\footnote{Amiens, 2 April 1869 - Paris, 29 April 1958, was a French physicist, engineer and mathematician known for his contributions to theoretical and applied physics, as well as to the theory of electrical circuits and dynamical systems}, who first studied it in the context of dynamical systems and differential equations. Li\'enard was known for his contributions in mathematical physics, particularly in the theory of electric circuits and nonlinear oscillations.

%\begin{figure}[h]
%\begin{center}
  %  \includegraphics{lienard.jpg}
  %\caption{Alfred-Marie Li\'enard}
%\end{center}
%\end{figure}

The Li\'enard equation is a second-order differential equation used to model nonlinear systems, and its general form is:

\begin{equation}\label{e:1}
\ddot x+f(x)\dot x+g(x)=0,
\end{equation}

where $\dot x$ represents the derivative of $x$ with respect to time, and $f(x)$ and $g(x)$ where $f$, $g$ are continuous functions $f,g:\mathbb{R}\rightarrow \mathbb{R}$ and $g(0)=0$.

This type of equation appears in the description of nonlinear oscillators in physics and electronics, such as in electric circuits that generate self-sustaining oscillations, such as relaxation circuits. Later, the Li\'enard equation gained further relevance when used in the analysis of the Van der Pol oscillator, a widely studied model in nonlinear dynamics and is the subject of detailed studies of many authors due to many applications in variety domains in science and technology, see for example \cite{GLN}, \cite{HHYS}, \cite{HZ}, \cite{H}, \cite{Las}, \cite{LNN}, \cite{N2}, \cite{N3}, \cite{PH} and classical sources \cite{Lie}, \cite{V1} and \cite{V2}.

Although Li\'enard did not develop all of the modern theory associated with his equation, his early work was fundamental to the analysis of stability and periodic oscillations in dynamical systems. The name thus pays tribute to his pioneering contribution to the subject.

The Li\'enard equation is important because it is a fundamental tool for analyzing nonlinear dynamical systems that exhibit self-sustaining oscillations, a phenomenon that appears in a wide range of scientific and technological disciplines. Here are some specific reasons that explain its relevance:

1. General Model of Nonlinear Oscillations 
The Li\'enard equation \eqref{e:1} generalizes many nonlinear systems with oscillatory behavior. It represents a broad class of physical and engineering problems, including electrical circuits, biological systems, and mechanical oscillations.

2. Applications in Electronic Circuits
It was key in the development of the theory of oscillators such as the Van der Pol oscillator, which is a special case of the Li\'enard equation.
It models phenomena such as oscillations in relaxation circuits and resonant systems.

3. Self-Sustaining Oscillations Phenomenon
The Li\'enard equation describes systems that, under certain conditions, exhibit limit cycles. These are stable periodic solutions that occur in nonlinear systems even in the absence of external forcing.
This concept is crucial in the study of natural oscillations in biology, such as cardiac or neuronal rhythms.

4. Stability Analysis
It provides a basis for studying the stability of periodic solutions and their dependence on initial conditions or system parameters.
It is fundamental to understanding the Hopf bifurcation, where limit cycles emerge from an equilibrium point.

5. Multidisciplinary Impact
In physics, it models mechanical systems with nonlinear friction or nonlinear resonance.
In biology, it appears in the modeling of metabolic cycles and neural signals.
In engineering, it is used in the design of oscillators for radio frequency and telecommunications systems.

6. Deep Mathematical Theory
It has driven the development of mathematical tools for the qualitative analysis of differential equations, such as phase diagrams, limit cycle theory, and the stability of nonlinear systems.
In summary, the Liénard equation not only describes fundamental oscillatory phenomena, but also provides a framework for studying nonlinear dynamical systems in a wide variety of fields. Its importance lies in its ability to model and analyze complex patterns that appear naturally in many scientific and technical contexts.

\section{Historic context}

Within the limits of Celestial Mechanics, Oscillation Theory and (partly) of mathematical physics, new integration methods were developed in the 20th century, apart from qualitative methods, in particular those that are solved with series according to special functions. This meant the development of the apparatus of Mathematical Analysis, thus completing the body of a mathematical theory related to ordinary differential equations and not, as had been happening until then, the development of the methods of ordinary differential equations, following the development of the mathematical theories themselves. This is an index, of course, of the quality of mathematical thought at the end of the last century.

On the other hand, in the study of certain physical systems, it is interesting, and almost always necessary, to know properties (of the solutions of the equation or system that models such system) such as boundedness, stability, periodicity, etc., without having to resort to the arduous and laborious task, which in many cases is impracticable, of finding analytical expressions for the solutions. In this way, the problem of investigating the properties of the solutions of a differential equation from ``its own expression” arose, giving rise to the Qualitative Theory of Differential Equations, a theory that emerged in the second half of the 19th century and was initially addressed by Jules Henri Poincar\'e (1854-1912) and Alexander Mijaílovich Liapunov (1857-1918) although for different reasons (one, due to the study of equilibrium figures and the stability of motion ``Probl\`eme g\'eneral de la stabilit\'e du mouvement”, French translation (1907) of the original in Russian (1893), ``Sur les figures d'equilibre peu diffèrentes des ellipsoids d'une masse liquide homogéne dovée d' un mouvement de rotation” (1906) and the other, due to his research in Celestial Mechanics, ``Sur les courbes définie par une équations differentielle” (1880), ``Memoire sur les courbes définie par une équation differentielle” (1881-1886) and ``Les Methodes Nouvelles de la Mecanique Celeste” (1892-1899).

The Second Lyapunov Method, also called the ``Direct Method”, is not only used nowadays to prove theorems of stability and theoretical mechanics. Currently, this method is applied to practical problems of mechanics and electrical oscillations, particularly, in engineering control. The theory of the direct method has received considerable progress in recent years and is approaching a certain state of completion. Some works have appeared in the last decades, related to the study of the behavior of the trajectories of nonlinear differential equations, which represent gradual generalizations of the Li\'enard equation \eqref{e:1}.

The method admits a very simple geometric interpretation, which is very likely due to G. N. Chetaev (1902-1959) and which is particularly useful in applications.

If we add to this that in the last 50 years we have witnessed a remarkable development of a field of Physics-Mathematics, designated by the name of Non-Linear Mechanics, a term, probably not entirely correct, since the changes have not occurred in Mechanics itself, but mostly in the techniques of solving its problems, especially those treated with the help of linear differential equations, which now use non-linear differential equations.

This is not a new idea in Mechanics. In fact, these non-linear problems are known from the studies of Euler, Poinsot, Lagrange and other geometers of that time, sufficient to illustrate this non-linear period of more than 150 years. The main difficulty of these studies, called classical, lies in the absence of a general method to treat these problems, which were treated above all, by special devices to obtain their solution.

It was Andronov (``Theory of Oscillations”, Moscow, 1929, written in collaboration with A. A. Witt and E. E. Chaikin\footnote{Alexander Vitt was arrested in 1937 during the Stalinist purges on charges of anti-Soviet activities. He died in custody shortly afterwards, and his name did not appear in the title of the first printed edition, the other authors citing a printing error. In later editions justice was done and his name appeared posthumously. A condensed English translation was published as part of the project on nonlinear differential equations under contract to the Office of Naval Research, edited under the direction of Solomon Lefschetz in 1947}) who suggested that self-oscillating oscillations are experienced in terms of the limit cycles (isolated periodic solutions) of Poincaré theory, marking a new era in these studies. Moreover, in 1937, he and L.S. Pontriagin (1908-1983) defined the notion of a rough (structurally stable) system, taking into account that one of the characteristics of the self-oscillating system may be the presence of the so-called feedback that controls the energy consumption of the non-periodic source, it follows directly from all the above that the mathematical model of such a system must be rough and notoriously non-linear. Intuitively, a dynamical system is structurally stable if small `perturbations' of the system do not change the phase state. Later, with the development of qualitative theory, it has been shown that structural stability is an important notion in this theory and has therefore been the subject of much research in recent years, in particular by the S. Smale school.

We implicitly distinguish four stages in the development of the Oscillation Theory, namely a first stage beginning with the work of Poincaré and ending in 1929; a second stage (developed almost exclusively in the former USSR with the work of Andronov and his school) spanning the period 1929-1940; the third stage, from 1945 and ending in the mid-1960s, and marked by the rapid development of the Theory of Automatic Control, so that in this stage a series of new problems appeared and the last stage, from the 1970s to the present day, which is greatly influenced by the development of electronic computing media, with all the advantages that this brings in the simulation and numerical solution of problems.

In this work, we will review some of the latest results obtained for various generations of the classical Li\'enard equation, using different generalized and fractional operators.

\section{The Li\'enard Equation}

In the strict sense of the word, the equation studied by Li\'enard was

\begin{equation*}
\ddot x+f(x)\dot x+x=0,
\end{equation*}

with the equivalent system

\begin{equation*}
\left\{ 
\begin{array}{ll}
\dot{x}=y-F(x), &  \\ 
\dot{y}=-x, & 
\end{array}
\right.
\end{equation*}

with $F'(x)=f(x)$. Li\'enard proved (see \cite{Lie}) that if F is a continuous odd function, which has a unique positive root at $x = a$ and is monotone increasing for $x \geq a$, then the above system,  has a unique limit cycle.

An additional detail, if we consider that $F$ is a polynomial, we are in the presence of a particular case of the second part of Hilbert's 16th problem: to find a bound for the number of isolated closed orbits for a polynomial vector field of degree $d$ on the plane. 

In these qualitative studies, one of the basic questions is the global existence, the prolongability of the solutions. For example, in \cite{N1} some properties of the system

\begin{equation}\label{e:2}
\left\{ 
\begin{array}{ll}
\dot{x}=\alpha (y)-\beta (y)F(x), &  \\ 
\dot{y}=-g(x), & 
\end{array}
\right.
\end{equation}

was studied under suitable assumptions. It is clear that under the conditions $\alpha(y)=y$, $\beta(y) \equiv 1$ of the system \eqref{e:2} we obtain the equation \eqref{e:1}.

Letting $F(x)=\int_{0}^{x}f(s)ds$ we obtain an equivalent system to equation \eqref{e:1}:

\begin{equation}\label{e:3}
\left\{ 
\begin{array}{ll}
\dot{x}=y-F(x), &  \\ 
\dot{y}=-g(x), & 
\end{array}
\right.
\end{equation}

Fractional calculus studies problems with derivatives and integrals of real or
complex order. As a purely mathematical field, the theory of fractional calculus was
brought up for the first time in the XVIIth century and since then many renowned
scientists worked on this topic, among them Euler, Laplace, Fourier, Abel, Liouville
and Riemann (see \cite{MKM}). Today there is a vast literature on this subject, many works and researchers multiply day by day in the world showing the most varied applications. The most common applications are currently found in Rheology, Quantum Biology, Electrochemistry, Dispersion Theory, Diffusion, Transport Theory, Probability and Statistics, Potential Theory, Elasticity, Viscosity, Automatic Control Theory, ...
 
After presenting the fractional integrals it is possible to define the fractional
derivatives associated with these operators. Fractional derivative and fractional integral are extensions of ordinary calculus, by considering derivatives of arbitrary real or complex order, and a general form for multiple integrals. Although mathematicians have wondered since the very beginning of calculus about these questions, only recently they have proven their
usefulness and since then important results have appeared not only in mathematics,
but also in physics, applied engineering, biology, etc. One question that is important is what type of fractional operator should be considered, since we have in hand several distinct definitions and the choice dependes on the considered problem.

In the literature many different types of fractional operators have been proposed, here,  we show that various of that  different notions of derivatives, can be considered particular cases of our definition and, even more relevant, that it is possible to establish a direct relationship between global (classical) and local derivatives, the latter not very accepted by the mathematical community, under two arguments: their local character and compliance with the Leibniz Rule.
To facilitate the understanding of the scope of our definition, we present the best known definitions of integral operators and their corresponding differential operators (for more details you can consult \cite {CT}). Without many difficulties, we can extend these definitions, for any higher order.

\section{A general integral operator}

We assume that the reader is familiar with the classic definition of the Riemann Integral, so we will not present it. In \cite{GLNV} was presented an integral operator generalized, which contain as particular cases, many of the well-known integral operators, both integer order and not. First we will present the definition of generalized derivative, see \cite{NGLK} (see also \cite {ZL} and \cite{FNRS}) which was defined in the following way.

\begin{definition} \label{d:1}
Given a function $f:[0,+\infty )\rightarrow \mathbb{R}$. Then the N-derivative of $f$ of order $\alpha $\ is defined by 

\begin{equation}\label{e:dn}
N_{F}^{\alpha }f(t)=\underset{\varepsilon \rightarrow 0}{\lim }\frac{
f(t+\varepsilon F(t,{\alpha }))-f(t)}{\varepsilon } 
\end{equation}

for all $t>0$, $\alpha\in (0,1)$ being $F(\alpha,t)$ is some function. 
If $f$ is $\alpha -$differentiable in some $(0,\alpha )$, and $\underset{t\rightarrow 0^{+}}{\lim }N_{F}^{\alpha}f(t)$ exists, then define $N_{F}^{\alpha}f(0)=\underset{t\rightarrow 0^{+}}{\lim }
N_{F}^{\alpha}f(t)$, note that if $f$ is differentiable, then $N_{F}^{\alpha}f(t)=F(t,\alpha)f'(t)$ where $f'(t)$ is the ordinary derivative.
\end{definition}

We consider the following examples:

I) $F(x,\alpha) \equiv 1$, in this case we have the ordinary derivative.

II) $F(x,\alpha )={ E }_{ 1,1 }({ x }^{ -\alpha  })$. In this case we obtain, from Definition \ref{d:1}, the non conformable derivative ${ N }_{ 1 }^{ \alpha  }f(x)$ defined in \cite{GLLMN} (see also \cite{NGL1}).

III) $F(x,\alpha)={ E }_{ 1,1 }((1-\alpha )x)={ e }^{ (1-\alpha)x }$, this kernel satisfies that $F(x,\alpha)\rightarrow 1$ as $\alpha \rightarrow 1$, a conformable derivative used in \cite{FMNS}.

IV)$F(x,\alpha)={ { E }_{ 1,1 }({ x }^{ 1-\alpha  }) }_{1}={ x }^{1-\alpha}$ with this kernel we have  $F(x,\alpha)\rightarrow 0$ as $\alpha \rightarrow 1$ (see \cite{KHYS}), a conformable derivative. 

V)$F(x,\alpha)={ { E }_{ 1,1 }({ x }^{ -\alpha  }) }_{1}={ x }^{\alpha}$ with this kernel we have  $F(x,\alpha)\rightarrow x$ as $\alpha \rightarrow 1$ (see \cite{NRS}). It is clear that since it is a non-conformable derivative, the results will differ from those obtained previously, which enhances the study of these cases.

VI)  $F(x,\alpha)={ { E }_{ 1,1 }({ x }^{ -\alpha  }) }_{1}={ x }^{-\alpha}$ with this kernel we have  $F(x,\alpha)\rightarrow x^{-1}$ as $\alpha \rightarrow 1$ This is the derivative ${ N }_{ 3 }^{ \alpha  }$ studied in \cite{MMN}. As in the previous case, the results obtained have not been reported in the literature.

Now, we give the definition of a general fractional integral. Throughout the work we will consider that the integral operator kernel $T$ defined below is an absolutely continuous function.

\begin{definition} \label{d:01}
Let $I$ be an interval $I \subseteq \RR$, $a,t \in I$ and $\a \in \RR$.
The integral operator $J_{F,a}^\a$, right and left, is defined for every locally integrable function $f$ on $I$ as

\begin{equation}\label {e:oig+}
J_{ F,a+ }^{ \alpha  }(f)(x)=\int _{ a }^{ x } \frac { f(t) }{ F\left( t-x,\alpha  \right)  } dt,\quad x>a.
\end{equation}

\begin{equation}\label {e:oig-}
J_{ F,b- }^{ \alpha  }(f)(x)=\int _{ x }^{ b } \frac { f(t) }{ F\left( x-t,\alpha  \right)  } dt,\quad x<b.
\end{equation}
\end{definition}

\begin{remark}
If not indicated, we will understand that the integral operator is defined over the entire interval, i.e., $J_{ F,a }^{ \alpha  }(f)(b)=\int _{ a }^{ b } \frac { f(t) }{ F\left( t,\alpha  \right)  } dt$.
\end{remark}

\begin{remark}
It is easy to see that the case of the $J_F^\alpha$ operator defined above contains, as particular cases, the integral operators obtained from classic, conformable and non-conformable local derivatives. However, Readers can see that many fractional integral operators can also be derived without much difficulty from the previous one.\newline
It is clear then, that from our definition, new extensions and generalizations of known integral operators can be defined. \newline
We can define the function space $L_{\alpha}^{ p}[a,b]$ as the set of functions over $[a,b]$ such that $({J}_{F,a+}^{\alpha}{[f(x)]^p}(b))<+\infty$.
\end{remark}

The following results are generalizations of the known results of the integer order Calculus.

\begin{proposition} \label{p:fundamental}
Let $I$ be an interval $I \subseteq \RR$, $a \in I$, $0<\alpha \le 1$ and $f$ a $\alpha$-differentiable function on $I$ such that $f'$ is a locally integrable function on $I$.
Then, we have for all $x \in I$
$$
J_{F,a+}^\a \big({ N }_{ F }^{ \alpha  }(f)\big)(x)=f(x)-f(a).
$$
\end{proposition}

\begin{proposition} \label{p:inverse}
Let $I$ be an interval $I \subseteq \RR$, $a \in I$ and $\a \in (0,1]$.
$$
N_F^{\a} \big( J_{F,a+}^\a(f)\big) (x)= f(x),
$$
for every continuous function $f$ on $I$ and $a,t\in I$.
\end{proposition}

In \cite{KHYS} it is defined the integral operator $J_{F,a}^\a$ for the choice of the function $F$ given by $F(x,\a)= x^{1-\a}$,
and \cite[Theorem 3.1]{KHYS} shows
$$
N^{\a} J_{x^{1-\a}\!,\,a}^\a(f)(x)= f(x),
$$
for every continuous function $f$ on $I$, $a,x \in I$ and $\a \in (0,1]$.
Hence, Proposition \ref{p:inverse} extends to any $F$ this important equality.

\begin{theorem} \label{t:pos}
Let $I$ be an interval $I \subseteq \RR$, $a,b \in I$ and $\a \in \RR$.
Suppose that $f,g$ are locally integrable functions on $I$, and $k_1,k_2\in \RR$. Then we have

$\eqref{e:1}$ $J_{F,a}^\a \big( k_1 f+k_2 g \big) (x) =k_1 J_{F,a}^\a f(x) + k_2 J_{F,a}^\a g(x),$

$\eqref{e:2}$ if $f \ge g$, then $J_{F,a}^\a f(x)\ge J_{F,a}^\a g(x)$ for every $t\in I$ with $t \ge a$,

$\eqref{e:3}$ $\left| J_{F,a}^\a f(x) \right| \le J_{F,a}^\a \left| f \right| (x)$ for every $t\in I$ with $t \ge a$,

$(4)$ $\int _{ a }^{ b } \frac { f(s) }{ T(s,\alpha ) } ds=J_{ T,a }^{ \alpha  }f(x)-J_{ T,b }^{ \alpha  }f(x)=J_{ T,a }^{ \alpha  }f(x)(b)$ for every $t \in I$.
\end{theorem}

Let $C^{1}[a,b]$ be the set of functions f with first ordinary derivative continuous on $[a,b]$, we consider the following norms on $C^{1}[a,b]$:

\begin{equation*}
{ \left\| F \right\|  }_{ C }=\max _{ [a,b] }{ \left| f(x) \right|  } ,\quad { \left\| F \right\|  }_{ { C }^{ 1 } }=\left\{ \max _{ [a,b] }{ \left| f(x) \right|  } +\max _{ [a,b] }{ \left| f'(x) \right|  }  \right\} 
\end{equation*}

\begin{theorem} (\cite{GLNV})
The fractional derivatives ${ N }_{ F,a+ }^{ \alpha  }f(x)$ and $ { N }_{ F,b- }^{ \alpha  }f(x)$  are bounded operators from $C^{1}[a,b]$ to $C[a,b]$ with

\begin{equation}
\left| { N }_{ F,a+ }^{ \alpha  }f(x) \right| \le K{ \left\| F \right\|  }_{ C }{ \left\| f \right\|  }_{ { C }^{ 1 } },
\end{equation}

\begin{equation}
\left| { N }_{ F,b- }^{ \alpha  }f(x) \right| \le K{ \left\| F \right\|  }_{ C }{ \left\| f \right\|  }_{ { C }^{ 1 } },
\end{equation}

where the constant K, may be depend of derivative frame considered.
\end{theorem}

\begin{remark}
From previou results we obtain that the derivatives ${ N }_{ F,a+ }^{ \alpha  }f(x)$ and $ { N }_{ F,b- }^{ \alpha  }f(x)$ are well defined.
\end{remark}

\begin{theorem} (Integration by parts) \label {t:parts}
Let $f,g:[a,b] \rightarrow \mathbb{R}$ differentiable functions and $\alpha \in (0,1]$. Then, the following property hold

\begin{equation}\label{e:8}
J_{ F,a+ }^{ \alpha  }((f)({ N }_{ F,a+ }^{ \alpha  }g(x)))={ \left[ f(x)g(x) \right]  }_{ a }^{ b }-J_{ F,a+ }^{ \alpha  }((g)({ N }_{ F,a+ }^{ \alpha  }f(x))).
\end{equation}
\end{theorem}

So, raises in natural way the study of more general system than \eqref{e:1}, within the framework of the N-derivative.

\section{Results}

\subsection{The non-conformable case}

We consider the following equation

\begin{equation}\label{e:ncl}
{ N }_{ e^{t^{-\alpha }} }^{ \alpha  }\left( { N }_{ e^{t^{-\alpha }} }^{ \alpha  }x \right) +f(x){ N }_{ e^{t^{-\alpha }} }^{ \alpha  }x+a(t)g(x)=0,
\end{equation}

where $a$, $f$ and $g$ are continuous functions satisfying the following conditions: \newline

a) $xg(x)\mathtt{>0}$ for $x\neq 0$, $g\in C^{1}(\mathbb{R}).$

b) $_{ { N }_{ e^{t^{-\alpha }} } }^{  }{ { J }_{ 0 }^{ \alpha  } }g(+\infty )=+\infty .$

c) $0<a\leq a(t)\leq A<+\infty $ for $t\in \left[ 0,+\infty \right) .\newline
$

The derivative ${ N }_{ e^{t^{-\alpha }} }^{ \alpha  }$ was defined ans studied in \cite{GLLMN} (see also \cite{NGL1}).

In \cite{NT} the follopwing results was obtained.

\begin{theorem}\label{t:o}
Let the following  conditions hold:
\begin{enumerate}
\item[1)] $_{ { N }_{ e^{t^{-\alpha }}} }{ { J }_{ 0 }^{ \alpha  } }\left( \frac { { N }_{ e^{t^{-\alpha }} }^{ \alpha  }a(t)^{ - } }{ a(t) }  \right) (+\infty )<\infty .$

\item[2)] $xg(x){ > }0,\quad x\neq 0.$

\item[3)] There is $N>0$ such that $\left\vert F(x)\right\vert \leq N$\ for $x\in\mathbb{R}$.
Then all solutions of the system (3) are oscillatory if and only if
\end{enumerate}
\begin{equation}\label{o:1}
_{ { N }_{ e^{t^{-\alpha }} } }{ { J }_{ { t }_{ 0 } }^{ \alpha  } }a(t)g\left[ \pm k(t-t_{ 0 }) \right] (+\infty )=\pm \infty ,
\end{equation}

for all $k\mathtt{>}0$ and all $t_{0}\geq 0$.
\end{theorem}

\begin{lemma} \label{l:1}
Under conditions a)-c) if $f(x)\in F_{g}(\mathbb{R})$ and ${ N }_{ e^{t^{-\alpha }} }^{ \alpha  }a(t)>0,$ then all solutions of system \eqref{e:ncl} are continuable to the future, i.e., for all $t\geq t_{0}\geq 0$.
\end{lemma}

\begin{theorem}\label{t:2}
Under assumptions of Lemma \ref{l:1} the following conditions:

\begin{enumerate}

\item[1)] ${ N }_{ e^{t^{-\alpha }} }^{ \alpha  }a(t)\mathtt{>}0$ for $t\geq 0$,
\item[2)]  $\left\vert F(x)\right\vert \leq N$ for some $N\mathtt{>}0$ and $x\in \mathbb{R}$,
\item[3)]  $G(\infty )=\infty ,$

\end{enumerate}

hold. Then the solutions of the equation \eqref{e:ncl} are bounded if and only if the
condition \eqref{o:1} is fulfilled.
\end{theorem}

\begin{lemma}\label{l:2}
If in addition to conditions of previous the theorem we have that g(x) is not
increasing function and $a(t)\rightarrow +\infty $ as $t\rightarrow +\infty $. Then condition \eqref{o:1} does not hold.
\end{lemma}

\begin{theorem}
Under condition Lemma \eqref{l:2} if the conditions

\begin{enumerate}
\item[a)] $_{ { N }_{ e^{t^{-\alpha }} } }{ { J }_{ 0 }^{ \alpha  } }\left( \frac { { N }_{ e^{t^{-\alpha }} }^{ \alpha  }a(t)^{ - } }{ a(t) }  \right) (+\infty )<\infty ,$
\item[b)] $G(x)\rightarrow +\infty $ as $\left\vert x\right\vert \rightarrow
+\infty ,$
\end{enumerate}

hold, then all solutions of equation \eqref{e:fle} are bounded.

\end{theorem}

\subsection{The generalized case}

In \cite{CN} the following sytems was studied:

\begin{equation}\label{e:gs}
\left. 
\begin{array}{l}
{\ N}_{F}^{\gamma }{x}=E(x,y) \\ 
{N}_{F}^{\gamma }{y}=-p(y)g(x)%
\end{array}%
\right\} \text{,}  
\end{equation}%
where $E(x,y)=a(x,y)H(x,y)$ with $H(x,y)=H(\alpha (y)-\beta (y)\Gamma(x))$ while the functions in \eqref{e:gs} are continuous real functions in their arguments.

\begin{remark}
The interested reader will be able to verify that, under different considerations on the functions involved in \eqref{e:gs}, the equations \eqref{e:1}, \eqref{e:2}, \eqref{e:3} and \eqref{e:fle} are obtained easily. Therefore, the results obtained generalize the previous ones.
\end{remark}

Throughout this paper we will use the following notation:

\begin{equation*}
\begin{array}{l}
G(x)=_{N}^{{}}{{J}_{F,0}^{\gamma }}(g)(x) \\ 
F(x)=_{N}^{{}}{{J}_{F,0}^{\gamma }}(f)(x) \\ 
E_{a,p}^{\alpha }(y)=E_{a,p}^{\gamma }(0,y)=_{N}^{{}}{{J}_{F,0}^{\gamma }}(%
\tilde{E})(y) \\ 
\tilde{E}(y)=\frac{a(0,y)H(\gamma (y))}{p(s)}%
\end{array}
\text{.}
\end{equation*}

Following the ideas in \cite{MN} the generalized partial derivatives can be
defined as follows.

\begin{definition}
Let a real valued function $f:\mathbb{R}^{n}\rightarrow \mathbb{\mathbb{R}}$
and $\overrightarrow{a}=(a_{1},\ldots ,a_{n})\in {\mathbb{R}}^{n}$ be point
whose $ith$-component is positive. Then the generalized partial N-derivative
of $f$ of order $\gamma $\ at the point $\overrightarrow{a}=(a_{1},\ldots
,a_{n})$ is defined by

\begin{equation}
N_{x_{i}}^{\gamma }f(\overrightarrow{a})=\underset{\varepsilon \rightarrow 0}%
{\lim }{\ \frac{f(a_{1},..,a_{i}+\varepsilon F(t,\gamma ),\ldots
,a_{n})-f(a_{1},...,a_{i},\ldots ,a_{n}))}{\varepsilon }}\text{.}  
\end{equation}

If that derivative exists, it is denoted by ${\ N}_{x_{i}}^{\alpha }f(\overrightarrow{a})$ and called the $ith$ non-conformable partial derivative of $f$ of order $\gamma \in (0,1]$ at $\overrightarrow{a}$.
\end{definition}

In this work, the following result were obtained, which generalize those of \cite{LNN}.

\begin{theorem} \label{t:1} 
Assume that:

i) $E(x,y)$ is a N-differentiable function on some finite open ball $B(x)$
around a point $x$ with 

\begin{equation*}
{\underset{B(x)}{sup}}{\ {\ N}_{x}^{\gamma }E(x,y)}<+\infty \text{ \ \ \ and
\ \ \ }\underset{y\rightarrow \pm \infty }{\lim \sup }\ E_{a,p}^{\alpha
}(y)=\pm \infty .
\end{equation*}

ii) $p(y)>0$ for all $y\in\mathbb{R}$.

iii) There exists a positive constant $\lambda $ such that 

\begin{equation*}
g(x)F(x)\geq -\lambda ,\text{ \ \ }_{N}{\ J}_{F,0}^{\gamma }{\ (xg)(t)}\geq
-\lambda \text{ for all }x\in \mathbb{R}\text{ \ \ \ and}\quad \underset{%
x\in \mathbb{R}}{\left\vert \Gamma(x)\right\vert }<+\infty .
\end{equation*}

iv) $a\in C^{1,0}(\mathbb{R}^{2})$ such that $a(x,y)>0$ for all $x$ and $y$
with

\begin{equation*}
\underset{x\rightarrow \pm \infty }{\lim \sup }\ a(x,y)<+\infty \text{ \ \ \
for all }y.
\end{equation*}

Then all solutions of system \eqref{e:gs} are defined for all $t$.
\end{theorem}

Based on this result and some complementary ones, several numerical simulations are presented for various $F$ kernels and different derivation orders.

\subsection{The fractional case}

In \cite{GLN2} the following system is studied:

\begin{equation}\label{e:fle}
\left\{ 
\begin{array}{ll}
_{C}^{0}\textrm{D}_t^\alpha {x}=y-F(x), &  \\ 
_{C}^{0}\textrm{D}_t^\alpha {y}=-g(x), & 
\end{array}
\right.
\end{equation}

where the functions involved satisfy the conditions considered above, and $_{C}^{0}\textrm{D}_t^\alpha $ is the Caputo-type Fractional Drift defined as follows:

\begin{definition} (Caputo fractional derivative).
The Caputo fractional derivative of order $\alpha $ on the half axis is
defined as follows

\begin{equation*}
_{a}^{C}D_{t}^{\alpha }f(t)=\frac{1}{\Gamma (n-\alpha )}
\int_{a}^{t} \frac{f^{(n)}(s)}{(t-s)^{\alpha +1-n}}ds, n-1<\alpha <n,  n\in \mathbb{N},\ 0<a<t<\infty .
\end{equation*}
\end{definition}

For any Lyapunov function $V\in C([t_{0},+\infty )$x$ \mathbb{R}^{n},\mathbb{R}_{+})$ we define

\begin{equation*}
 ^{C}D_{+}^{q}V(t,x)=\underset{h\rightarrow 0}{\lim \ \sup }\ \frac{1%
}{h^{q}}\left[ V(t,x)-V(t-h,x-h^{q}f(t,x)\right]
\end{equation*}

for $(t,x)\in [t_{0},+\infty )\times \mathbb{R}^{n}$. 

After presenting some definitions for the system

\begin{equation}
_{0}^{C}D_{t}^{\alpha }x(t)=f(t,x(t))
\end{equation}

Various results are obtained for this, in particular the following Lemma, basic for the subsequent results:

\begin{lemma}
Let $x(t)$ be a continuous and derivable function, $F(x)=\int_{0}^{x}f(r)dr $ and $f:\mathbb{R}
\rightarrow \mathbb{R}$ be a positive function satisfying a Lipschitz condition $\left\vert f(x(t))-f(x(\tau ))\right\vert \leq M\left\vert x(t)-x(\tau )\right\vert$ for some constant $M>0$.
Then, for any time instant $t\geq 0$ we have

\begin{equation*}
_{0}^{C}D_{t}^{\alpha }F(x(t))\leq f(x(t))_{0}^{C}D_{t}^{\alpha }x(t),\alpha \in (0,1).
\end{equation*}
\end{lemma}

Which are applied to the fractional Li\'enard type equation \eqref{e:fle}, so we have:

\begin{theorem}
Let $x(t)$ and $y(t)$ continuous and derivable functions. If $f$ and $g$ are
continuous functions satisfying:

i) $f:\mathbb{R}\rightarrow \mathbb{R}_{+}$ a Lipschitzian function,

ii) $g:\mathbb{R}\rightarrow \mathbb{R}$ and $xg(x)>0$ for $x\neq 0$,

then, the trivial solution of system \eqref{e:fle} is stable.
\end{theorem}

\section{Conclusions}

In this review we have presented some of the results obtained in recent years, concerning the qualitative behavior of solutions of the Li\'enard Equation, which generalize the known classical results.

We have highlighted these generalizations for the case of conformable, generalized and fractional operators of the Caputo type.

All this can be used as an abbreviated compendium, which can save time and effort to interested researchers.

\end{document}